# GLOBAL SUM ON SYMMETRIC NETWORKS

*Vance Faber*

**Abstract.** *We discuss the problem global sum on computer networks. Each processor starts with a single real value. At each time step, every directed edge in the network graph can simultaneously be used to transmit a single number between the processors (vertices). How many time steps $s$ are required to ensure that every processor acquires the global sum?*

Our general model of a network is a directed graph with processors as vertices and the connections between them as edges. The case has been made elsewhere that a multiprocessor network should be homogeneous; that is, the network should appear the same from any processor. This means that the graph is vertex transitive. Sabiduissi [1] has shown that a graph is vertex transitive if and only if it is the Cayley coset graph of a group. If the coset is the identity group, the graph is called a Cayley graph.

We are interested in the following problem we call *global sum*. Each processor starts with a single real value. At each time step, every directed edge in the graph can simultaneously be used to transmit a single number between the processors (vertices). The processors retain any information that they obtain during this process and there is ostensibly no restriction on what functions the processors perform to generate the number that they will transfer.

**Problem A.** *How many time steps $s$ are required to ensure that every processor acquires the global sum?*

This problem has similarity to the gossip problem (see [2]) but the underlying graph is not complete and we allow use of all the connections simultaneously. It is easy to see that a lower bound on the number of time steps needed to perform the task is the diameter $D$ of the graph. In addition, we do not intend to allow the length of information exchanged between the processors on a single time step to grow longer and longer as is the case in the ordinary gossip problem. Without this limitation, it is easy to see that the exact number of steps required is $D$. On every time step, each processor accepts all the information that its neighbors have and concatenates it together. After $D$ steps, everyone has all the information. For example, the number that is exchanged between processors could consist of blocks of bits, one block assigned to each processor. As soon as a processor discovers one of the numbers belonging to a processor, it could write that into the appropriate location. In order to forbid that, we assume that the number of bits allowed to be transferred between processors at any time is roughly on the order of the number of bits required to store the product of the number of processors times the size of the largest allowed value.

This problem is quite general. For example, if the values are allowed to be vectors, then by weighting the values assigned to the processors we might be asking for the number of communication steps required for each processor to acquire a matrix times a vector.

We have discussed this problem in [3]. There we mentioned that an upper bound is clearly $2D$. We span the graph with a tree of diameter $D$ with root $v$. The values are then communicated along the tree inward to $v$ while adding at each branch point. Then the global sum is broadcast from $v$ back to all the other vertices (possibly on a different tree). We optimistically conjecture that the true value is $D$ for all vertex symmetric graphs. We shall call a graph *sum optimal* if $s = D$. We showed in [3; Theorem 5.1] that the degree of the minimal polynomial of the adjacency matrix can be used to provide a different upper bound.

**Theorem 1.** *Suppose that the regular graph of degree $d$ with adjacency matrix $A$ has exactly $m+1$ distinct eigenvalues. Then the time for global sum is at most $m$.*

**Proof.** We utilize the Hoffman polynomial (see [4; page 157]). We choose time steps of the form $(A - \lambda_t I)$ where the $\lambda_t$ are the eigenvalues other than $d$ (which corresponds to the constant eigenvector $e = \vec{1}/\sqrt{n}$). After the $k^{th}$ time step, the $i^{th}$ processor has the $i^{th}$ entry of the vector $\prod_{t=1}^{k}(A - \lambda_t I)x$. Now if we expand $x$ into its parts along $e$ and orthogonal to $e$, we have $x = y + (x,1)\vec{1}/n$ where the sum of the entries of $y$ is zero. Furthermore $\prod_{t=1}^{D}(A - \lambda_t I)y = 0$, so after $m$ time steps, the entries are all the same value

$$\mu = \frac{(x,1)}{n} \prod_{t=1}^{m}(d - \lambda_t).$$

We recover the global sum at each processor by scaling $\mu$ by the predetermined factor.

Using this theorem, we also showed that weakly distance regular graphs satisfy the conjecture. Weakly distance regular graphs are a generalization of distance regular graphs and include all distance transitive graphs.

**Corollary 2.** *If all the edges in a weakly distance regular graph can be used simultaneously, then global sum takes $D$ time steps where $D$ is the diameter.*

**Proof.** A weakly distance regular graph has exactly $D+1$ eigenvalues (see [5, page 113]).

**Theorem 3.** *Suppose that $G_1$ and $G_2$ are two sum optimal Cayley graphs with underlying groups $\Gamma_1$ and $\Gamma_2$ and generating sets $\Delta_1$ and $\Delta_2$. Then the Cayley graph $G$ on the group $\Gamma_1 \times \Gamma_2$ with generating set $\Delta_1 \times \{e\} \cup \{e\} \times \Delta_2$ is sum optimal.*

**Proof.** Clearly, the distance between two vertices in $G$, $(a,b)$ and $(c,d)$ is the sum of the distances between $a$ and $c$ and $b$ and $d$ so the diameter of $G$ is just the sum of the diameters of $G_1$ and $G_2$. If we run the optimal summation algorithm on all the copies of $G_1$ followed by running the optimal summation algorithm on all the copies of $G_2$, each vertex will have the sum of all the values.

This theorem can be used to show that there exist vertex symmetric graphs that are sum optimal but not distance regular. Let $G_1$ be the 5-cycle and $G_2$ be a single edge. Then the product $G$ is not distance regular (see [6]).

Discussion. Consider general algorithms that can be expressed by restricting the functions performed by the processors to be linear combinations of the incoming values. Let $J$ be the order $n$ matrix all of whose entries is 1. In particular, since asking for each processor to acquire the global sum is asking for each processor to compute an equation in the system $y = Jx$, we ask this question.

**Problem B.** *What is the smallest $m$ so that $\prod_{t=1}^{m} W_t = J$, where each $W_t$ is a matrix with non-zero values only where the adjacency matrix is non-zero?*

Discussion. The number of unknowns is $mdn$ while the number of equations is $n^2$ so in any graph where $D < n/d$ it is not obvious why there would be a solution with $m = D$. However, one could say something similar for $m = 2D$ yet we know there is a solution in that case so these $n^2$ equations are not completely independent.

In the case where the graph is the Cayley graph of a group $\Gamma$, we can ask for the matrices $W_t$ to have the extra property that the nonzero entries $W_t(g, h)$ depend only on the generator $gh^{-1} = \delta$. If in addition, the group is a cyclic group, then the matrices are *circulant* matrices. In this case, each $W_t$ has the same eigenvectors, namely the discrete Fourier vectors $f_j$. Since the eigenvalues of $J$ are all zero except for the first which corresponds to the constant vector $f_0$, for each other eigenvector $f_j$ there must exist a $t$ such that $W_t f_j = 0$. Using the convolution theorem, we can then convert this problem into a problem solely about vectors.

**Problem C.** Let $V$ be the integers modulo $n$. Let $S$ be a subset of $d+1$ integers in $V$ including 0. What is the smallest collection of non-zero vectors $w_1, w_2, ..., w_m$ indexed by $V$ and with support $S$ so that for each non-constant Fourier vector $f_j$ there exists a $w_k$ with the inner product $(w_k, f_j) = 0$?

The following was shown in [7; Theorem 3.2].

**Theorem 4.** *Let A be the adjacency matrix of a graph G. Let $p_m(x)$ be a polynomial of degree m with $p_m(0) = 1$ such that*

$$\left\| p_m(A-dI) - \frac{J}{n} \right\| < \frac{1}{n-1}.$$

*Then G has a diameter $D \leq m$.*

We are not able to prove that this condition also implies that $s \leq m$. However, we can show that after $s$ time steps, each processor can calculate an approximate value for the global sum.

**Theorem 5.** *Let A be the adjacency matrix of a graph G. Suppose each vertex $i$ is assigned a value $x_i$. Let $p_m(x)$ be a polynomial of degree m with $p_m(0)=1$ such that*

$$\left\| p_m(A-dI) - \frac{J}{n} \right\| < \frac{1}{n-1}.$$

*Then after m communication steps, the vertex $i$ can calculate an approximate value $y_i$ for the mean $\mu$ of the $x_i$ with the property that the relative variance is less than $1/(n-1)^2$.*

**Proof.** Let $y = p_m(A-dI)x$. Then since $\frac{J}{n}x = \mu$, we have

$$\|y - \mu\| < \frac{1}{n-1}\|x\|.$$

**Theorem 6.** *Suppose each vertex $i$ is assigned a value $x_i$. Then after 2 communication steps, each vertex $i$ can calculate the sum of all the $x_j$ of distance 2 from $i$.*

**Proof.** On the first time step, each vertex $i$ exchanges its value with its neighbors. Let $j$ be adjacent to $i$ and $k$ be adjacent from $i$ but not from $j$. Let $n(j,k)$ be the number of paths of length 2 from $j$ to $k$. On time step 2, vertex $i$ transmits the value

$$\sum_{(j,i)} \frac{x_j}{n(j,k)}$$

to vertex $k$. Then vertex $k$ adds all the values it received on time step 2. At this point, the value that vertex $k$ has is given by

$$\sum_{(i,k)}\sum_{(j,i)} \frac{x_j}{n(j,k)} = \sum_{k} \sum_{\substack{i \\ d(j,k)=2 \ (j,i),(i,k)}} \frac{x_j}{n(j,k)} = \sum_{\substack{k \\ d(j,k)=2}} \frac{x_j}{n(j,k)} \sum_{\substack{i \\ (j,i),(i,k)}} 1 = \sum_{\substack{k \\ d(j,k)=2}} x_j.$$

**Corollary 7**. *A graph of diameter 2 is sum optimal*.

**Proof**. Each vertex needs only to add its own value and that of its neighbors on the last step of the previous algorithm.

**Remark**. I don't see how to extend this to diameter 3. The difference may lie in the fact that this algorithm cannot be expressed in the terms of Problem B.